%%%%%%%%%%%%%%%%%%%%%%%%%%%%%%%%%%%%%%%%%%%%%%%%%%%%%%%%%%%%%%%%%%%%%%%%%%
%%%% Just a TRANSLATION with some nonessential changes from:
%%%%   Reinov O.I., Approximation of operators in Banach spaces,
%%%%   Primenenie funkcional'nogo analiza v teorii priblizhenij (1985), p. 128-142
%%%%   (in Russian).
%%%%%%%%%%%%%%%%%%%%%%%%%%%%%%%%%%%%%%%%%%%%%%%%%%%%%%%%%%%%%%%%%%%%%
%%%%%   Note:
%%%%% We describe the topology $ \tau_p$ in the space $ \Pi_p(Y,X),$
%%%%% for which the closures of convex sets in $\tau_p$ and in
%%%%% $ {}^*$-weak topology of the space $ \Pi_p(Y,X)$ are coincident.
%%%%% Thereafter, we investigate some properties of the space $ \Pi_p,$
%%%%% related to this new topology. 2010-remark:
%%%%% Occasionally, the topology is coincides with the $\lambda_p$-topology
%%%%% from the paper "Compact operators which factor through subspaces of $l_p$",
%%%%% Math. Nachr.  281(2008), 412-423) by Deba Prasad Sinha and Anil Kumar Karn.
%%%%
%%%%%%%%%%%%%%
%
% This is pdfTeX, Version 3.1415926-1.40.10 (MiKTeX 2.8)
% (preloaded format=amstex 2010.1.26)
% 14 FEB 2010 18:20
% AmS-TeX- Version 2.2
%
%%%%%%%%%%%%%%%%%%%%%%%%%%%%%%%%%%%%%%%%%%%%%%%%%%%%%%%%%%%
%
%%%%%%%%
%%%%%%%%%%%%%%%%%%%%%%%%%%%%%%%%%%%%%%%%%%%%%%%%%%%%%%%%%%%%%%%%
%%%%    Kalinin 1985 Primenen. FA v t. pribl.
%%%%           31.01.00 16:31:28 Sun
%%%%%%%%%%%%%%%%%%%%%%%%%%%%%%%%%%%%%%%%%%%%%%%%%%%%%%%%%%%%%%%%
%%%%%%%%%%%%%%
 \input amstex
\documentstyle{amsppt}

 \expandafter\edef\csname aaaaa \endcsname{%
  \catcode`\noexpand\@=\the\catcode`\@\space}
\catcode`\@=11
\def\nologo{\let\logo@\empty}
\csname aaaaa \endcsname

\nologo
% \Monograph
   %\input trtrtrtr.tex
\magnification=\magstep1
\parindent=1em %%%% ?
         %\vsize=7.4in   %%%% 1in = 2.54
 %\NoPageNumbers
     \baselineskip=18pt      %%%%%

\hsize=16 true cm
\vsize=24 true cm
%\mathsurround=3pt
                  %\pagewidth{11.5 true cm}
                  %\pageheight{17.00 true cm}
                  %\tenpoint

 %
    %
           %
\def\ove#1{\overline{#1}}    %
\def\ovs#1#2{\overset{#1}\to{#2}}
     \def\({\left(}       \def\al{\alpha}           \def\lee{\leqslant}
     \def\){\right)}      \def\e{\varepsilon}    \def\gee{\geqslant}
     \def\[{\left[}       
     \def\]{\right]}      \def\ffi{\varphi}

                                      \def\ot{\otimes}
     \def\<{\langle}                 \def\wh{\widehat}
     \def\>{\rangle}                 \def\wt{\widetilde}
                 \def\sbs{\subset}
\def\tr{\operatorname{trace}\,}

  %%%%%%%%%%%%% after 30.01.00 02:44:40 Sat:

\def\AP{\operatorname{AP}}
\def\BAP{\operatorname{BAP}}
\def\N{\operatorname{N}}
\def\I{\operatorname{I}}
\def\id{\operatorname{id}}
\def\L{\operatorname{L}}
\def\QN{\operatorname{QN}}

\def\R{\operatorname{R}}

                     \def\sbs{\subset}

 %%%%%%%%%%%%%%%%%%%%%%%%%

 %       \def\{\quad\blacksquare}

    \def\QQ{$\quad\blacksquare$}      
    \def\Q{\quad\blacksquare}         \def\bigp{\bigpagebreak}

\CenteredTagsOnSplits                        \NoBlackBoxes

%\,

%\vskip1in

\NoRunningHeads
\pageno=1  %
\footline={\hss\tenrm\folio\hss} %
    %
    % \ifnum\pageno<0 \romannumeral-\pageno \else\number\pageno \fi

        \topmatter
        \title {Approximation of operators in Banach spaces}
        \endtitle
          \author {O.I. Reinov}  \endauthor
\address\newline
                          %\tenrm
198904, Saint Petersburg,\newline
St. Petersburg State University,\newline
Dept. Math.,\newline
%Chair of Math. Analysis
\endaddress

\email
                          \tenrm
orein51\@mail.ru
\endemail

\abstract\nofrills
\endabstract
    \endtopmatter
   %==============================
 \footnote""{
[R--KGU]\ Reinov O.I., Approximation of operators in Banach spaces,
in book Primenenie funkcional'nogo analiza v teorii priblizheni,
Kalinin: KGU, 1985, 128-142.
  }

\document

 \head \S 0. Notations\endhead

This work is the first part of some investigations which are concerned
with an approximation of ones or others classes of operators in Banach
spaces (the approximation can be understood in different sences).
The main ojects of the study will be the spaces, possessing (or not
possessing) the so-called approximation perperties $\AP_p$ and
the bounded approximation properties $\BAP_p$ of order $p,$ where
$p>0$ (see, e.g., [9], [10]).
In the propounded part, we consider only the properties $\AP_p$
where $p\gee 1,$ and our considerations are bounded by some equivalent
reformulations of the first-given definition from [9].
The main accent is put on the investigation of the conditions, under which
the finite dimensional operators are dense in norm in the space of
quasi-$p$-nuclear operators, as  well as in the space of all
absolutely $p$-summing operators
in the topology of "$\pi_p$-compact convergence".

We will keep the standard notations and the terminology.
If $A$ is a bounded subset of a Banach space $X$, then
$\ove{\Gamma(A)}$ is the closed absolutely convex hull of $A;$
$ X_A$ is the Banach space with "the unit ball" $ \ove{ \Gamma(A)}$;
$ \Phi_A: X_A\to X$ is the canonical embedding. For $ B\sbs X,$
it is denoted by $ \overline{B}^{\,\tau}$ and $ \ove{B}^{\,\|\cdot\|}$
the closures of the set $ B$ in the topology $ \tau$ and in the norm $ \|\cdot\|$
respectively. When it is necessary, we denote by $ \|\cdot\|_X$
the norm in $X.$ Other notations:\,
$ \Pi_p,$ $ \QN_p,$ $ \N_p,$ $ \I_p$ are the ideals of absolutely $p$-summing,
quasi-$p$-nuclear, $p$-nuclear , strictly $p$-integral operators,
respectively (see [8], [10]);
$ X^*\wh\ot_p Y$ is the complete tensor product associated with
$ \N_p(X,Y);$
$ X^*\wh{\wh\ot}_p Y= \ove{X^* \ot Y}^{\,\pi_p}$ \, (i.e. the closure of the set
of finite dimensional operators in $ \Pi_p(X,Y)).$ Finally,
if $ p\in[1,+\infty]$ then $ p'$ is the adjoint exponent.

Everywhere bellow, we use to suppose that $ p\in [1,+\infty];$ however, the proofs
are given for $ p\in (1,+\infty)$\, (with some non-essential changes, all proofs
pass through the cases $ p=1$ and $ p=+\infty).$
The domain of changes for $p$ is noted specially only in the places where it is
necessary, or where a vagueness can be arised.

%%%%%%%%%%%%%%%%%%%%%%%%%%%%%%%%%%%%%%%%%%%%%%%%%%%%%%%%%%%%%%%
%\head
%  \bf Aproximation Properties:\, Topological Aspect. With Appendix
%    to a Grothendieck Conjecture
%\endhead
%\med
%
%
%\bigp
%%%%%%%%%%%%%%%%%%%%%%%%%%%%%%%%%%%%%%%%%%%%%%%%%%%%%%%%%%

 \head \S 1. Approximation of absolutely $p$-summing operators\endhead

In this paragraph we will describe the topology $ \tau_p$ in the space $ \Pi_p(Y,X),$
 for which the closures of convex sets in $\tau_p$ and in
 $ {}^*$\!-weak topology of the space $ \Pi_p(Y,X)$ are coincident
(this is an answer to a corresponding question of P. Saphar in [12], p.385).
Thereafter, we will investigate some properties of the space $ \Pi_p,$
related to this new topology.

Thus, consider in the space $ \Pi_p(Y,X)$ the {\it topology $ \tau_p$\ of
$ \pi_p$\!-compact convergence}, a local base (in zero) of which
is defined by sets of type
$$ \omega_{K,\e}= \left\{ U\in \Pi_p(Y,X):\ \pi_p(U\Phi_K)<\e\right\},
$$
where $ \e>0,$\, $ K=\ove{\Gamma(K)}$ --- a compact subset of $ Y.$

\proclaim {\bf 1.1. Proposition.}\it               %Teorema 1 !!!
Let $ \R$ be a linear subspace in
$ \Pi_p(Y,X),$  containing
$ Y^*\ot X.$ Then
$ (\R,\tau_p)'$ is isomorphic to a factor space of the space
$ X^*\wh\ot_{p'} Y.$
More precisely, if $ \ffi\in (\R, \tau_p)',$
then there exists an element $ z=\sum_1^\infty x'_n\ot y_n\in X^*\wh\ot_{p'} Y$
such that
$$ \ffi(U)= \tr\, U\circ z,\ \, U\in \R. \tag$*$
$$
On the other hand, for every
$ z\in X^*\wh\ot_{p'} Y$ the relation
$ (*)$ defines a linear continuous functional on
$ (\R,\tau_p).$
\endproclaim\rm

          %%%%%%%%%%%%%%%%%%%%%%%%%%%%%%%%%%%%%%%%%%%% OK till here !!!
The\, {\it Proof}\ of this proposition will be given in detail.
However, we will omit details in analogues cases hereinafter.
So, let
$ \ffi$ be a linear continuous functional on $ (\R,\tau_p).$
Then one can find a neighborhood of zero
$ \omega_{K}=\omega_{K,\e},$ such that $ \ffi$ is bounded on it:
$ \forall\, U\in\omega_K, \ |\ffi(U)|\lee 1.$
We may assume that $ \e=1.$ Consider the operator
$ U\Phi_K:\,Y_K\ovs{\Phi_K}\longrightarrow  Y
               \ovs{U}\longrightarrow X.$
Since the mapping  $ \Phi_K$ is compact,  $ U\Phi_K\in \QN_p(Y_K,X).$
Put $ \ffi_k(U\Phi_K)=\ffi(U)$ for $ U\in \R.$
On the linear subspace
$ \R_K= \left\{ V\in \QN_p(Y_K,X):\ V=U\Phi_K \right\}$
of the space  $ \QN_p(Y_K,X),$
the linear functional
$ \ffi_K$
is bounded: if
$ V=U\Phi_K\in \R_K$ and $ \pi_p(V)\lee 1,$
then $ |\ffi_K(V)|=\ffi(U)|\lee 1.$
Therefore, $ \ffi_K$ can be extended to a linear continuous functional
$ \wt\ffi$ on the whole  $ \QN_p(Y_K,X);$ moreover, because of the injectivity of the ideal
$ \QN_p,$
considering $ X$ as a subspace of some space $ C(K),$
we may assume that $ \wt\ffi\in\QN_p(Y_K, C(K))^*.$ Let us mention that
$$ \wt\ffi(jU\Phi_K)=\ffi_K(U\Phi_K)=\ffi(U) \tag1
$$
(here $ j$ is an isometric embedding  of $ X$ into $ C(K)$).

Furthermore, since
$ \QN_p(Y_K, C(K))^*= \I_{p'}(C(K), (Y_K)^{**}),$
we can find an operator
$ \Psi: C(K)\to (Y_K)^{**},$
for which
$$ \wt\ffi(A)= \tr \Psi A,\ \, A\in (Y_K)^*\ot C(K).
$$

Let
$ A_n\in (Y_K)^*\ot C(K),\, $ $ \pi_p(A_n-jU\Phi_K)\to 0.$
Then
$$ \wt\ffi(jU\Phi_K)= \lim\, \tr \Psi A_n. \tag2
$$

Consider the operator
$ \Phi_K^{**}\Psi: C(K)\ovs{\Psi}\longrightarrow (Y_K)^{**}
   \ovs{\Phi_K^{**}}\longrightarrow Y.$
Since
$ \Psi\in \I_{p'},$ and $ \Phi_K$
is compact, we have
$ \Phi_K^{**}\Psi \in \N_{p'}(C(K), Y)= C(K)^*\wh\ot_{p'} Y.$
% Ho!: go via 2 compact operators and refl. space by Johnson!!!
Let
$ \sum_1^\infty \mu_n\ot y_n\in C(K)^*\wh\ot_{p'} Y$
be a representation of the operator
$ \Phi^{**}_K\Psi.$
Put
$ z=\sum j^*(\mu_n)\ot y_n.$
The element
$ z$
generates an operator
$ \Phi_K^{**}\Psi j$
from
$ X$
to
$ Y.$
We will show now that
$ \tr U\circ z=\wt\ffi(jU\Phi_K)$
(note that
$ U\circ z$
is an element of the space
$ X^*\wh\ot X,$
so the trace is well defined).
We have:
$$\multline
 \tr U\circ z= \tr \( \sum j^*(\mu_n)\ot Uy_n\)=
   \sum \< j^*(\mu_n), Uy_n\> =\\ =
\sum \< \mu_n, jUy_n\>=
\tr jU\Phi_K^{**}\Psi=\tr (jU\Phi_K)^{**}\Psi,
\endmultline                                     \tag3
$$
where
$ (jU\Phi_K)^{**}\Psi:
    C(K)\ovs{\Psi}\longrightarrow (Y_K)^{**}
         \ovs{\Phi_K^{**}}\longrightarrow Y
       \ovs{U}\longrightarrow X\ovs{j}\longrightarrow C(K).$
% identify X with its image in  X**; come through with a compact act!!!!
Since
$ \pi_p(A_n- jU\Phi_K)\to 0,$
then
$ \pi_p\(A^{**}_n - (jU\Phi_K)^{**}\)\to 0.$
Moreover, if
$A:= A_n=\sum_1^N w_m\ot f_m\in (Y_K)^*\ot C(K),$
then
$$ \tr A_n^{**}\Psi= \sum_m \< \Psi^*w_m, f_m\> =
   \sum_m \< w_m, \Psi f_m\>= \tr \Psi A.
$$

Hence,
$ \tr (jU\Phi_K)^{**}\Psi = \lim\,\tr A^{**}_n\Psi=
\lim\,\tr \Psi A_n.$
Now, it follows from (3) and (2) that
$ \wt\ffi(jU\Phi_K)=\tr U\circ z.$
Finally, we get from (1):
$ \ffi(U)= \tr U\circ z.$
Thus, the functional
$ \ffi$
is defined by an element of
$ X^*\wh\ot_{p'} Y.$

Inversely, if
$ z\in X^*\wh\ot_{p'} Y,$
put
$ \ffi(U)=\tr U\circ z$
for
$ U\in \R$
(the trace is defined since
$ U\circ z\in X^*\wh\ot X).$
We have to show that the linear functional
$ \ffi$
is bounded on a neighborhood
$ \omega_{K,\e}$
of zero in
$ \tau_p.$
For this, we need

\proclaim {\bf 1.2. Lemma}\it
If
$ z\in X^*\wh\ot_q Y,$
then $ z\in X^*\wh\ot_q Y_K,$
where
$ K=\ove{\Gamma(K)}$
is a compact in
$ Y.$
\endproclaim\rm

\demo{\it Proof of the lemma
}\
Let
$ z=\sum x'_n\ot y_n,$ $ \{c_n\}\in c_0$
and
$ \sum \|x'_n\|^q\,c_n^{-1}<+\infty. $
Consider the operator
$ A_1:l_q\to Y,$ $ A_1\{a_n\}= \sum a_n c_n y_n.$
Since this operator is compact, one can find a compact
$ K\sbs Y$
and an operator
$ A_2: l_q\to Y_K,$
for which
$ A_1=\Phi_K A_2.$
Put
$ z_0=\sum c_n^{-1}\,x'_n\ot e_n$\, ($ e_n$
are orths in
$ l_q).$
Then
$ \ove{z}:=(\bold 1\ot A_2)(z_0)\in X^*\wh\ot_q Y_K,$
and
$ \Phi_K(\ove{z})=z.$
$\quad\blacksquare$\enddemo   %!!!

Let us continue the proof of the theorem. Let
$ K$ be a compact subset of
$ Y,$
for which
$ z\in X^*\wh\ot_{p'} Y_K.$
If
$ U\in \omega_{K,1},$
then
$ \pi_p(U\Phi_K)<1$
and
$ |\tr U\Phi_K\circ z|\lee
\|z\|_{X^*\wh\ot_{p'} Y_K}\cdot$  %\,
$\pi_p(U\Phi_K)\lee C.$
$\quad\blacksquare$

%%%%%%%%%%%%%%%%%%%%%%%%%%%%%%%%

\proclaim {\bf 1.3. Corollary}\it
$ (\R,\tau_p)'=(\R,\sigma)',$
where
$ \sigma=\sigma(\R, X^*\wh\ot_{p'} Y).$
Thus, the closures
of convex subsets of the space
$ \Pi_p(Y, X)$
in
$ \tau_p$
and in
$ \sigma$
are the same.
$\quad\blacksquare$  \endproclaim\rm

Denote by
$ X^*\wt\ot_{p'} Y$
the closure of the set
$ X^*\ot Y$
in the space
$ \I_{p'}(X, Y^{**})$
(dual to
$ Y^*\wh{\wh\ot}_p X).$

\proclaim {\bf 1.4. Proposition}\it
Let
$ A$
be the intersection of the unit ball of the (dual to
$ X^*\wt\ot_{p'} Y)$
space
$ G=G(Y, X^{**})$
with the subspace
$ Y^*\ot X.$
$ {*}$\!-weak
closure of the set
$ A$ in $ G\cap \Pi_p(Y,X)$
coincides with the closure of
$ A$
in
$ (\Pi_p(Y, X), \tau_p).$
\endproclaim\rm

\demo{\it
Proof
}\
Let us consider the canonical mappings
$$ \CD
  X^*\wh\ot_{p'} Y  @>j>>  X^*\wt\ot_{p'} Y, \\
\Pi_p(Y, X^{**})    @<j^*<<  G(Y, X^{**}).
\endCD
$$

Since
$ j^*$
is one-to-one, then the closures of the bounded sets in
$ G(Y, X^{**}),$
in topologies
$ \sigma(G, X^*\wt\ot_{p'} Y)$
and
$ \sigma (G, X^*\wh\ot_{p'} Y)$
are the same. Therefore, if
$ B$ is a convex bounded set in
$ G(Y,X)\sbs \Pi_p(Y,X),$
then the closure of the set
$ B$
in
$ \( \Pi_p(Y,X)\cap G, \sigma(G, X^*\wt\ot_{p'} Y)\)$
coincides with the closure of the set
$ B$
in the space
$ \( \Pi_p(Y,X), \sigma(\Pi_p(Y,X), X^*\wh\ot_{p'} Y)\)$
and, therefore,
by Corollary 1,3, --- with the closure of
$ B$
in
$ \( \Pi_p(Y,X), \tau_p\).$
\QQ \enddemo

                    %%%%%%%%%%%%%%%%%%%%%%%%%%%%%%% OK till now 12:40 PM 11/1/2009 ???

\proclaim {\bf 1.5. Corollary}\it
With notations of the proposition {\rm 1.4},
the closure of the set
$ A$
in
$ \tau_p$
coincides with the closure of
$ A$
in the space
$ \L(Y,X)$
in the topology of compact convergence.
\endproclaim\rm

For the\, {\it proof,}\ it is enough to use the previous assertion, considering
the canonical mapping from
$ X^*\wh\ot_1 Y$
into
$ X^*\wt\ot_{p'} Y$
instead of the map
$ j$
from the proof of the proposition {\rm 1.4}
(and to apply either Proposition {\rm 1.1}
for
$ p=+\infty,$
or results on duality from [4]).
$\quad\blacksquare$

%%%%%%%%%%%%%%%%%%%%%%%%%%%%%%%%%%%%%%%%%%%%%%%%%%%!!!
%%%%%%%%%%%%%%%%%%%%%%%%%%%%%%%%%%%%%%%%%%%%%%%%%%%!!!
\proclaim {\bf 1.6. Corollary}\it
Let
$ C>0$
and
$ T\in \Pi_p(Y,X).$
The following assertions are equivalent:

$1)$
there ia a net
$ \left\{ T_\al\right\}, T_\al\in Y^*\ot X,$
converging to
$ T$
in the topology
$ \tau_p$
such that
$ \pi_p(T_\al)\lee C;$

$2)$
there is a net
$ \left\{ T_\al\right\}, T_\al\in Y^*\ot X,$
converging to
$ T$
in the topology of compact convergence,
such that
$ \pi_p(T_\al)\lee C.$
$\quad\blacksquare$
\endproclaim\rm

\proclaim {\bf 1.7. Proposition}\it
For an operator
$ T\in \Pi_p(Y,X), \ove{T(Y)}=X,$
the following are the same:

$1)$ $ T\in \ove{Y^*\ot X}^{\,\tau_p};$

$2)$
there is a net of operators
$ R_\al\in Y^*\ot Y$
such that
$ TR_\al\to T$
in the topology
$ \tau_p.$
\endproclaim\rm

\demo{\it Proof}\
Assuming that 2) is not valid, we (by Corollary 1.3) get:
$$ \not\exists \, R_\al\in Y^*\ot Y:\ TR_\al\to T\ \text{ in }
   \( \Pi_p(Y,X), \sigma(\Pi_p(Y,X), X^*\wh\ot_{p'} Y)\).  \tag4
$$

Consider the associated with
$ T$
mappings:
$$ \CD
  X^*\wh\ot_{p'} Y @>\wt{T}>>   Y^*\wh\ot_1 Y, \\
   \Pi_p(Y, X^{**})  @<\wt{T}^{\,*}<<  \L(Y, Y^{**}).
\endCD
$$
where
$ \wt{T}(z)= z\circ T$
for
$ z\in X^*\wh\ot_{p'} Y.$
Let
$ Z=\ove{\wt{T}^{\,*}(Y^*\ot Y)}^{\,*}$
(the closure is taken in the space
$ \Pi_p(Y,X)$
in
$ {}^*$\!-weak topology).
It follows from (4) that
$ T$
is not zero on the subspace
$ Z^{\perp}\sbs X^*\wh\ot_{p'} Y,$
i.e. there exists an element
$ A\in Z^{\perp}$
such that
$ \< T,A\>=\tr AT=1.$  But $ AT=\wt{T}(A),$
and if
$ R\in Y^*\ot Y,$ then $ \< \wt{T}(A), R\>=
\< A, (\wt{T})^{\,*}(R)\>=0.$
Hence, the element
$ \wt{T}(A)$
of the projective tensor product
$ Y^*\wh\ot Y$
is not zero
(since
$ \tr \wt{T}(A)=1),$
but generates a null-operator in $ Y.$
For any
$ y'\in Y^*$
and
$ Ty\in T(Y)$
we have:
$ \< A, Ty\ot y'\>= \< ATy,y'\>=0.$
Since
$ \ove{T(Y)}=X,$
we obtain that a non-zero tensor element
$ A\in X^*\wh\ot_{p'} Y$
generates a zero-operator. Again, by using the equality
$ \tr AT=1,$
we conclude that
$ T$
can not be approximated in
$ {^*}$\!-weak
topology by finite rank operators.
Now, it follows from Corollary 1.3 that %the condition
1)
is not fulfilled.
$\quad\blacksquare$
\enddemo

%%%%%%%%%%%%%%%%%%%%%%%%%%%%%%%%%%%%%%%  !!! 5-6

Next two statements give us sufficient (but not necessary, as we will see below)
conditions for the density of the set of all finite rank operators in the space
of operators
$ \Pi_p(Y,X)$
in the topology
$ \tau_p$\ of
$ \pi_p$-compact convergence.

\proclaim {\bf 1.8. Proposition{\rm\footnote{
In the original paper [R--KGU] this Proposition was sounded
as follows:\newline
{\phantom{AAA} }\qquad  "{\it If\ $ \QN_p(Y,X)= \ove{Y^*\ot X}^{\,\pi_p},$
then\  %!!!
   $ \Pi_p(Y,X)= \ove{Y^*\ot X}^{\,\tau_p}.$}"
      }}}
\it
%If $ \QN_p(Y,X)= \ove{Y^*\ot X}^{\,\pi_p},$ then  !!!
%$ \Pi_p(Y,X)= \ove{Y^*\ot X}^{\,\tau_p}.$
% --- So was formulation in that paper what is problematic -
%  so, in the case $p=\infty$ it is obtained somewhat open question
%    Here is a more normal formulation:           !!!
If
$ \QN_p(Y,X)= \ove{Y^*\ot X}^{\,\pi_p}$
for every Banach space
$ Y,$  then for each
$ Y$\
$ \Pi_p(Y,X)= \ove{Y^*\ot X}^{\,\tau_p}.$
\endproclaim\rm

\demo{\it Proof}
Let
$ U\in \Pi_p(Y,X),$ $ \e>0,$ $ K=\ove{\Gamma(K)}$
be a compact in
$ Y.$
By the assumptions, there is an operator
$ V\in (Y_K)^*\ot X,$
such that
$ \pi_p(V-U\Phi_K)<\e.$
      %  (this follows from
      %                                    [Theorem I.4.1.1)].    %%%%% ???
%  To clear this - or more general formulation on product of
%  summing and RN-operators!!!
We need to set successfully instead of
$ V$
an operator
$ \wt V\Phi_K,$
where
$ \wt V\in Y^*\ot X.$
Let
$ V=\sum_{n=1}^N z_n\ot x_n.$
Note that we can consider only the case where
$ \Phi_K^{**}$
is one-to-one
(else, with the help of the construction of [2]
we change
$ Y_K$
by a space
$ Y_{K_0},$
for which the operator
$ \Phi_{K_0}$
is compact and the operator
$ \Phi_{K_0}^{**}$
is one-to-one).
In this case
$ Y^*$  is norm dense in
$ (Y_K)^*$
and, therefore, for every positive number sequence
$ \{ \e_n\}$
there exist the elements
$ y'_n\in Y^*,$
for which
$ \|y'_n-z_n\|_{Y^*_K}<\e_n.$
Put
$ \wt V=\sum_1^N y'_n\ot x_n\in Y^*\ot X.$
Let
$ \{ a_n\}$
be a sequence of the elements of the space
$ Y_K,$
such that
$ \sup \left\{ \sum | \< a_n,a'\>|^p:\ \|a'\|_{Y^*_K}\lee 1\right\}\lee 1.$
We have:
$$ \multline
 \sum_{i=1}^m \|(\wt V\Phi_K-V)a_i\|^p =
  \sum_{i=1}^m \|\sum_{n=1}^N \< y'_n-z_n, a_i\>\,x_n\|^p\lee \\ \lee
  \sum_{i=1}^m \sum_{n=1}^N | \< y'_n-z_n, a_i\>|^p
    \( \sum_{n=1}^N \|x_n\|^{p'}\)^{p/{p'}}\lee  \\ \lee
    \( \sum_{n=1}^N \|x_n\|^{p'}\)^{p/{p'}}
   \sum_{n=1}^N \|y'_n-z_n\|^p_{Y^*_K}\lee
     \( \sum_{n=1}^N \|x_n\|^{p'}\)^{p/{p'}}
      \sum_{n=1}^N \e_n^p.
\endmultline
$$

If we take
$ \e_n$
small enough then the last number is less then
$ \e,$
and, from the inequality
$ \pi_p(V-U\Phi_K)\lee\e,$
we get that
$ \pi_p(\wt V\Phi_K-U\Phi_K)\lee \e+ \pi_p(V-\wt V\Phi_K)\lee 2\e.$
Hence,
$ \wt V-U\in \omega_{K,2\e}.$
Thus, we have shown that for every neighborhood
$ \omega_{K,\e}$\
there exists an operator
$ \wt V\in Y^*\ot X,$
for which
$ \wt V-U\in \omega_{K,\e}.$
Therefore
$ U\in \ove{Y^*\ot X}^{\,\tau_p}.$
$\quad\blacksquare$\enddemo

\proclaim {\bf 1.9. Proposition}\it
If the canonical mapping
$ j: X^*\wh\ot_{p'} Y\to \N_{p'}(X,Y)$
is one-to-one then
$ \Pi_p(Y,X)= \ove{Y^*\ot X}^{\,\tau_p}.$
\endproclaim\rm
\demo{\it Proof}
If the map
$ j$
is one-to-one then the annihilator
$ j^{-1}(0)^\perp$
of its kernel in the space, dual to
$ X^*\wh\ot_{p'}Y,$
coincides with
$ \Pi_p(Y, X^{**}).$
On the other hand, in any case
$ j^{-1}(0)^\perp=
\ove{Y^*\ot X}^{\,*}$
(the closure in
$ {}^*$\!-weak topology of the space
$ \Pi_p(Y, X^{**}));$
by Corollary 1.3,
$$ \Pi_p(Y,X)\cap \ove{Y^*\ot X}^{\,*}= \ove{Y^*\ot X}^{\,\tau_p}.
$$

Therefore,
$ \Pi_p(Y,X)= \ove{Y^*\ot X}^{\,\tau_p}.$
$\quad\blacksquare$\enddemo

For a reflexive space
$ X,$
the dual space to
$ X^*\wh\ot_{p'} Y$
is equal to
$ \Pi_p(Y,X).$
Consequently, it follows from 1.3 and 1.9

\proclaim {\bf 1.10. Corollary}\it
For a reflexive space
$ X$
the canonical mapping
$ j: X^*\wh\ot_{p'} Y\to \N_{p'}(X,Y)$
is one-to-one iff the set of finite rank operators is dense in the space
$ \Pi_p(Y,X)$
in the topology
$ \tau_p$\ of $ \pi_p$-compact convergence.
\QQ\endproclaim\rm

To conclude this part of our considerations, let us give an assertion
which shows the following. It  follows from the existing of
an absolutely $p$-summing operator, which is non-approximated in the topology
$ \tau_p,$
that there esits a non-approximated (in the same topology)
quasi-$p$-nuclear operator (with values in the same space).

\proclaim {\bf 1.11. Proposition}\it
If there exists an operator
$ T\in \Pi_p(Y,X)\setminus \ove{Y^*\ot X}^{\,\tau_p},$
then there exist a reflexive space
$ Z$
and an operator
$ U\in\QN_p(Z,X),$
which is not in the closure
$\ove{Y^*\ot X}^{\,\tau_p}.$
\endproclaim\rm
For the\,  {\it proof}\  it is enough to remember the definition of the topology $ \tau_p$
and to use the following two facts:\
a)\, if
$ V$
is a compact operator then it can be represented as a composition of two compact operators;\
b)\, the product
$ AB$
of a compact operator
$B$
and absolutely $ p$-summing operator
$ A$
is a quasi-$p$-nuclear map
(see [5], [8]).
\QQ

\remark {\bf 1.12. Remark}\footnote{In the original paper [R--KGU]
this non-essential remark is sounded exactly as
"At the moment, I do not know whether the inverse of Proposition 1.8
is true in the case where the space $X$ is reflexive.
  }     %%%%%%%%%%%  see Rus Tex and comment there !!! ???
\endremark\medpagebreak

     %%%%%%%%%%%%%%%%%%%%%%%%%%%%%%%%%%%%%%%%%%%%%%%%%%% OK ALL !!! ???? 1:47 PM 11/1/2009

         %----\S2 in App_85
\head \S 2. Approximation properties of Banach spaces \endhead

In this paragraph, we fix the space of images of operators
(or the spaces where operators are defined),
and investigate the conditions under which it is possible
to approximate (by finite rank operators) {\it all}\ absolutely $p$-summing mappings
with values in a given space (or acting from a given space).\footnote{
  In the original paper [R--KGU], there was here the phrase\
   "The next statement, among other things, gives us a partial inversion of
   Proposition 1.8."   %%% Same remark as just above in Rem 1.12.!!!        "
           }

\proclaim {\bf 2.1. Proposition}\it
For a Banach space
$ X$
the following are equivalent:

$1)$
for every Banach space
$ Y$
the equality
$ \QN_p(Y,X)= \ove{Y^*\ot X}^{\,\pi_p}$
holds;

$2)$
for every Banach space
$ Y$
the equality
$ \Pi_p(Y,X)= \ove{Y^*\ot X}^{\,\tau_p}$
holds;

$3)$
for every Banach space
$ Y$
one has %the equality
$ \QN_p(Y,X)\subset \ove{Y^*\ot X}^{\,\tau_p};$\footnote{
  In the original paper [R--KGU], there was here an evident misprint:
   "for every Banach space $ Y$ the equality
    $ \QN_p(Y,X)= \ove{Y^*\ot X}^{\,\tau_p}$ holds;"
          }

$4)$
for every reflexive Banach space
$ Y$
the equality
$ \QN_p(Y,X)= \ove{Y^*\ot X}^{\,\tau_p}$
holds.
\endproclaim\rm

\demo{\it Proof}
Implications
$ 2)\implies 3)\implies 4)$
are evident;
$ 1)\implies 2)$
by Proposition 1.8.
For the proof of the implication
$ 4)\implies 1),$
consider an operator
$ V\in \QN_p(Y,X).$
By using the results of
[2],   % [66],
 we can factorize the operator
$ V$
by the following way:
$ V=UA,$
where
$ A\in \operatorname{L}_c(Y,Z),$
         %(compact),
$ U\in \QN_p(Z,X),$
and, moreover, the space
$ Z$ is reflexive. Put
$ K=\ove{A(\operatorname{ Ball}_Y)}.$
From 4), it follows that
$$ \forall\,\e>0\ \exists\, \wt V\in Z^*\ot X:\ \ \wt V-U\in\omega_{K,\e},
$$
i.e.
$ \pi_p(\wt V\Phi_K-U\Phi_K)<\e.$
Hence,
$ \pi_p(V-\wt VA)\lee C_A\,\pi_p(U\Phi_K- \wt V\Phi_K)\lee C_A\e.$
$\quad\blacksquare$ \enddemo

\proclaim {\bf 2.2. Corollary}\it
If for every reflexive space
$ Y$
the canonical mapping
$ X^*\wh\ot_{p'} Y \to \N_{p'}(X,Y)$
is one-to-one then for each Banach space
$ Y$
the equality
$ \QN_p(Y,X)= \ove{Y^*\ot X}^{\,\pi_p}$
holds.
$\quad\blacksquare$ \endproclaim\rm

For the  {\it proof,}\,
 Proposition 1.9 may be applied, and then use
the implication
$ 2)\implies 1)$
of the previous fact.
 \QQ

It follows from Corollaries 1.10 and 2.2

\proclaim {\bf 2.3. Corollary}\it
For a reflexive Banach space
$ X$
the following are equivalent:

$1)$
for each space
$ Y$
the canonical mapping
$ X^*\wh\ot_{p'} Y \to \N_{p'}(X,Y)$
is one-to-one;

$2)$
for each space
$ Y$
the set of finite rank operators is dense in the Banach space
$ \QN_p(Y,X).$
$\quad\blacksquare$ \endproclaim\rm

The next statement, among other things, gives us a partial inversion of
Proposition 1.9.

\proclaim {\bf 2.4. Corollary}\it
For every Banach space
$ X$
the following are equivalent:

$1)$
for each space
$ Y$
the canonical mapping
$Y^*\wh\ot_p X \to \N_p(Y,X)$
is one-to-one;

$2)$
for each reflexive Banach space
$ Y$
the canonical mapping
$Y^*\wh\ot_p X \to \N_p(Y,X)$
is one-to-one;

$3)$
for each (reflexive) Banach space
$ Y$\
$ \ove{X^*\ot_{p'} Y}^{\,\tau_{p'}}= \Pi_{p'}(X,Y).$
\endproclaim\rm
\demo{\it Proof}\,
Concerning the equivalence
$ 1)\iff 2),$
see [1];  % [57]
 %or
 %Theorem 1;
the implications
$ 1)\implies 3)$
and
$ 3)\implies 2)$
follow from 1.9 and 1.10, respectively.
$\quad\blacksquare$ \enddemo

The previous statement yields the following well known result:

\proclaim {\bf 2.5. Corollary}\it
If a Banach space
$ X$
has the approximation property then for any
$ p>1$
and any space
$ Y$
the canonical mapping
$Y^*\wh\ot_p X \to \N_p(Y,X)$
is one-to-one.
$\quad\blacksquare$ \endproclaim\rm

    %%%%%%%%%%%%%%%%%%%%%%%%%%%%%%%%%%%%%%%%%%%%%%%%%%%% !!! OKnow ????? 3:04 PM 11/1/2009

The above results lead us to the following definition which is equivalent to
corresponding definition in [9], [10], [11].                 %%%%%%%%%%%%

\definition {\bf 2.6. Definition}
Let
$ p\gee 1.$
A Banach space
$ X$
has the property
$ \AP_p$
(the approximation property of order
$ p$),     %%%%%%%%%
if every absolutely
$ p'$-summing operator, acting from the space
$ X,$
can be approximated in the topology of
$ \pi_p'$-compact convergence
$ \tau_{p'}$
by operators of finite rank.
\enddefinition

It follows from  Proposition 1.8 and Corollary 2.4
 that

\proclaim {\bf 2.7. Corollary}\it
If for each Banach space
$ Y$
the equality
$ \QN_{p'}(X,Y)=\ove{X^*\ot Y}^{\,\pi_{p'}}$
holds then the space
$ X$
has the property
$ \AP_p.$
       %% in the sense of Definition $1.$
$\quad\blacksquare$ \endproclaim\rm

Recall for the sake of completeness the following
assertion on a characterization of the spaces with the property
$ \AP_p$
(a proof can be found in
[1]).        %%

\proclaim {\bf 2.8. Proposition}\it
A Banach space
$ X$
has the property
$ \AP_p$
   %% in the sense of Definition $1$ !!!
iff
for every Banach space
$ Y,$
for every operator
$ T\in \Pi_{p'}(X,Y),$
for each weakly
$ p'$-summable sequence
$ \{ x_k\}$
of elements of the space
$ X$
and for every
$ \e>0$
there is a finite rank operator
$ R: X\to Y,$
such that
$ \sum \|Ux_k-Rx_k\|^{p'}<\e.$
\endproclaim\rm

   %%%%%%%%%%%%%%%%%%%%%%%%%%%% ti;; Now 3:43 PM 11/1/2009 All OK but later... ??????!!!

As the property $ \AP_1$ (the usual approximation property og Grothendieck),
the properties $ \AP_p$
are very useful in investigations of the questions of different kinds
in the geometrical theory of operatos
(for instance, when describing the dual spaces of some spaces of operators;
so, if $ X$ possesses the property $ \AP_p,$ then
$ \N_p(Y,X)^*=\Pi_{p'}(X,Y^{**})$ for every space $ Y).$
However, we would like (concluding this paragrap) to adduce
a simple fact, which is valid without any assumptions on approximation
(see [3] for $p=1$ and [12] for $p>1,$ where an analogues assertion was proved
with a supposition of approximation property; see also [7],
Theorem 4.6, where it is obtained a little bit less general result).

\proclaim {\bf 2.9. Proposition}\it
Let $ p\in[1,+\infty).$
The following are equivalent

$1)$ Banach spaces $ X$ and $ Y$ are reflexive;

$2)$ the space $ \QN_p(X,Y)$ is reflexive;

$3)$ the space $ \Pi_p(X,Y)$ is reflexive.
\endproclaim\rm

\demo{\it Proof}\
Since, for reflexive spaces $ X$ and $ Y,$ the equality
$ \QN_p(X,Y)=\Pi_p(X,Y)$ holds, it is sufficient to prove only
the implication $ 1)\implies 2).$
For this, imbedd the space $ Y$ into some space $ C(K)$ isometrically,
and let us consider the space $ \QN_p(X,Y)$ as a subspace of
$ \QN_p(X,C(K)).$
Then any continuous functional $ \Phi$ on $ \QN_p(X,Y)$
has an extension with the same norm to a linear functional
$ \wt\Phi$ on $ \QN_p(X, C(K)).$ In turn, the functional $ \wt \Phi$
is generated by an operator $ U\in \I_{p'}(C(K),X^{**}),$ so
$ \< \Phi, T\>= \tr UjT$ for $ T\in \QN_p(X,Y)$ (where $ j$
is the imbedding of
$ Y$ into $ C(K)$). Since the space $ Y$ is reflexive, one has
$ Uj\in \N_{p'}(Y,X)$ [7]. Therefore, as a functional,
$ \Phi$ is generated by an element of the space $ Y^*\wh\ot_{p'} X.$
Thus, the natural mapping
$ Y^*\wh\ot_{p'} X \to \QN_p(X,Y)^*$ is an "onto" map.
Since  $ \(Y^*\wh\ot_{p'} X \)^*=\QN_p(X,Y),$
we get that the space $ \QN_p(X,Y)$ is reflexive.
$\quad\blacksquare$ \enddemo

                  %%%%%%%

  \head \S 3. Counterexamples \endhead

A lot of counterexamples concerning AP's
can be found in [9], [10], [11].
We will use them partially. Recall that
for any
$ p\gee 1,$ $ p\neq2,$
there exists a separable reflexive Banach space without the property
$ \AP_p.$
Moreover, for every
$ p\gee 1,$ $ p\neq2,$
there exist a separable reflexive
$ E,$
an operator
$ R\in \Pi_{p'}(E,E)$
and a tensor element
$ t\in E^*\wh\ot_p E,$
so that
$ \tr t\circ R=1$ and $ \tr t\circ A=0$
for each finite rank operator
$ A\in E^*\ot E$\, (for details, we refer the reader to the papers
[9], [10], [11]).

It is often very useful to apply the following fact when constructing
some counterexamples (see [6]):   %

\proclaim {\bf 3.1. Lemma}\it
For every separable Banach space $ E$ there exist a separable conjugate
Banach space $ H=Y^*$ with a basis and operators
$ Q: H\to E$ and $ U: H^*\to E^*$ so that $ Q(H)=E,$
$ U(H^*)=E^*,$ $ \|U\|\lee 1,$ $ \|Q\|\lee 1,$ $ UQ^*=\id_{E^*}$
and the space $ \( \id_{E^*}- Q^*U\)(H^*)$
is isomorphic to the space $ Y.$
   %%%
\endproclaim\rm

    %%%%%%%%%%%%%%%%%%%%%%%%%%%%%%%
   %%%%%%%%%%%%%%%%%%%%%%%%%%%%%%%%%%%%%%%%%%%%%

Now we are ready for constructions of our counterexamples.
Firstly, we will show that the inversion of
Proposition 1.9 and Corollary 2.2
are not valid.

\proclaim {\bf 3.2. Proposition}\it
For every
$ p\in [1,+\infty],$ $ p\neq2,$
there exist a separable reflexive space
$ E,$
a separable conjugate space
$ H$
with a basis such that the canonical mapping
$ j: H^*\wh\ot_{p'} E\to \N_{p'}(H,E)$
is not one-to-one. On the other hand, since
$ H$
has the AP,
$ \QN_p(E,H)=\ove{E^*\ot H}^{\,\pi_p}$
and
$ \Pi_p(E,H)= \ove{E^*\ot H}^{\,\tau_p}.$
\endproclaim\rm

\demo{\it Proof}
Let
$ E,$ $ t\in E^*\wh\ot_{p'} E$
and $ R\in \Pi_p(E,E)$
be the mentioned above spaces, tensor element and operator so that
$ \tr t\circ R=1$
and
$ t=0$
as an operator.
Set
$ g=t\circ Q\in H^*\wh\ot_{p'} E$
and
$ V=U^*R\in \QN_p(E, H^{**}).$
Then
$ \tr V\circ g=1$
(consequently, the tensor element
$ g$
is not equal to zero)
and
$ g=0$
as an operator.
$\quad\blacksquare$ \enddemo

\proclaim {\bf 3.3. Corollary}\it
For every
$ p\in[1,+\infty],$ $ p\neq2,$
there exist a Banach space
$ Z,$
an element
$ z\in Z^*\wh\ot_{p'} Z$
and an operator
$ \Psi\in\Pi_p(Z,Z^{**})$
such that
$ \tr \Psi\circ z=1,$
but
$ \tr \Phi\circ z=0$
for all
$ \Phi\in \Pi_p(Z,Z).$
\endproclaim\rm
\demo{\it Proof}\,
Let us use notation introduced in the proof of
Proposition 3.2.
Let
$ \sum y'_n\ot y_n$
be any representation of a tensor element
$ g\in H^*\wh\ot_{p'} E.$
Note that
$ \tr A\circ g=0$
for every operator
$ A\in \Pi_p(E,H)$
(because of the space
$ H$
has approximation property).
Put
$ Z=H\oplus E,$
$ z=\sum (y'_n, 0)\ot (0,y_n)$
and define the operator
$ \Psi\in\QN_p(Z,Z^{**})$
by
$ \Psi(h,y)=(Vy,0).$
We have:
$$ \tr \Psi\circ z=\sum \< (y'_n,0), (Vy_n, 0)\>=
     \sum \< y'_n, Vy_n\> =\tr V\circ g=1.
$$

On the other hand, denoting by
$ P_H$
and
$ P_E$
the natural projectors from
$ Z$
onto
$ H$
and
$ E$
respectively, we have, for arbitrary operator
$ \Phi\in \Pi_p(Z,Z),$

$$ \Phi(h,y)=\Phi|_H(h)+\Phi|_E(y)=
   \( P_H\Phi|_H(h)+P|_E\Phi|_H(h)\)+
   \( P_H\Phi|_E(y)+P|_E\Phi|_E(y)\),
$$
whence,
$$ \tr \Phi\circ z= \sum \< (y'_n, 0), \[ P_H\Phi(0, y_n)\]\>=
    \sum \< (y'_n,0), P_H\Phi|_E(y_n)\>.
$$
Denoting by
$ A$
the operator
$ P_H\Phi|_E,$
we get:
$$ A\in \Pi_p(E,H);\ \ \tr \Phi\circ z=\tr A\circ g=0. \Q
$$
\enddemo

\remark {\bf 3.4. Remark}
For
$ p=+\infty$
we get nonzero tensor element
$ z\in Z^*\wh\ot_1 Z,$
vanishing on the subspace
$ \L(Z,Z)$
of the space
$ \L(Z,Z^{**}).$
This is an answer to a question of Swedish mathematician
Sten Kaijser, who was one who is directed my attention
for a possibility of the
existemce of such an element $z.$
\endremark\medpagebreak

Now we will show that the inversion of
Proposition 1.8\footnote{
In [R--KGU] I meant that, for fixed $X$ and $Y,$\
in Proposition 1.8\newline
\phantom{AAAAa}$ \Pi_p(Y,X)= \ove{Y^*\ot X}^{\,\tau_p}\ \nRightarrow\
 \QN_p(Y,X)= \ove{Y^*\ot X}^{\,\pi_p}
$
                }
and Corollary 2.7 are invalid too.

\proclaim {\bf 3.5. Proposition}\it
For every
$ p\in [1,+\infty],$ $ p\neq2,$
there exist a separable reflexive space
$ E,$
a separable conjugate space
$ H$
with a basis
$($so, with the property
$ \AP_{p'})$
such that
$ \QN_{p}(H,E)\neq \ove{H^*\ot E}^{\,\pi_p};$
on the other hand,
$ \Pi_p(H,E)=\ove{H^*\ot E}^{\,\tau_p}$ \, {\rm (see Corollary 2.4).}
\endproclaim\rm
\demo{\it Proof}
With the notation of the proof of Proposition 3.2,
set
$ L=RQ.$
Since
$ \tr (t\circ RQ^{**}Y^*)=1$
and $ t=0$
as an operator,
the map
$ RQ^{**}$
can not be approximated by finite rank operators in the space
$ \QN_p(H^{**},E).$
Moreover,
$ L\not\subset \ove{H^*\ot E}^{\,\pi_p}.$
$\quad\blacksquare$
\enddemo

%  %%%%%%%%%%%

In conclusion, let us bring the following, at first sight somewhat surprising,
statement which shows, roughly speaking, that there are spaces $ X $
and $ Y, $ for which the closures in $\(\Pi_p(X,Y^{**}), w^*\)$
of the set of all finite dimensional operators is minimal: coincides with
$ X^*\wh{\wh\ot}_p Y. $

\proclaim {\bf 3.6. Proposition} \it
For every $ p\in [1,+\infty)$ there exist (separable and reflexive)
Banach spaces $ X$ and $ Y$
such that
$$ \ove{X^*\ot Y}^{\,\tau_p}= X^*\wh{\wh\ot}_p Y.
$$
 \endproclaim\rm
  \demo {\it Proof}\
 Let $ X $ and $ Y $ be th separable and reflexive spaces such
that the natural mapping
$ Y^*\wh\ot_{p'} X \to \N_{p'}(Y,X)$ is
not one-to-one. By Proposition 1.1, the dual space to
$ \N_{p'}(Y,X)$ can be identified with a subspace
$ \ove {X ^*\ot Y} ^ {\, \tau_p} $
of the spaces $ \Pi_p (X, Y) $
(see also the proof of Proposition 1.9).
Since this subspace is reflexive
(Proposition 2.9) and
$ (X^*\wh{\wh\ot}_p Y)^*=\N_{p'}(Y,X),$
we have:
$X^*\wh{\wh\ot}_p Y = \ove{X^*\ot Y}^{\,\tau_p}.$
 $ \quad\blacksquare $
\enddemo

\remark {\bf 3.7. Remark} Seemingly, it is unknown
whether Proposition 3.6 is true in the case where $ p=+\infty.$
    \endremark

\bigp

%%%%%%%%%%%%%%

%%%%%%%%%%%%%%%%%%%%%%%%%%%%%%%%%%%%%%%%%%%%%%%%%%%%%%%%%%%%%

\Refs


\ref \no1 \by Bourgain J., Reinov O.I. \pages 19-27
\paper On the approximation properties   for the space $H^\infty$
\yr 1985\vol   122
\jour Math. Nachr.
\endref


\ref \no2 \by  Davis W.J, Figiel T., Johnson W.B., Pelczynski A.
                             \pages     311-327
\paper    Factoring weakly compact operators
\yr 1974 \vol 17 \issue
\jour J. Functional Analysis
\finalinfo     $MR 50 \# 8010.$
\endref

\ref \no 3  \by  Y.Gordon, D.R.Lewis, H.R.Retherford \pages 85-129
\paper   Banach ideals of operators with applications
\yr 1973 \vol 14 \issue 1
\jour      J. Funct. Anal.
\finalinfo
\endref



\ref \no4 \by Grothendieck A. \pages pp. 196 + 140
\paper  Produits tensoriels topologiques et espases nucl\'eaires
\yr 1955\vol  16
\jour  Mem. Amer. Math. Soc.
\finalinfo
\endref


\ref \no5
\by Johnson W.B. \pages 337-345
\paper  Factoring compact operators
\yr 1971 \vol 9 \issue
\jour         Israel J. Math.
\endref


\ref \no6 \by Lindenstrauss J.\pages  279-284
\paper  On James' paper "Separable Conjugate Spaces"
\yr 1971\vol 9
\jour Israel J. Math.
\endref

 \ref \no 7\by Makarov B.M., Samarskij V.G.
 \paper  Weak sequencial completeness and close properties of
he spaces of operators
 \inbook  Theory of operators and theory of functions
\publaddr Leningrad
\publ  LGU
 \yr 1983\vol\nofrills 1,
\pages\nofrills   122-144
 \endref


\ref \no8
\by   Pietsch A.
\book  Operator ideals
\bookinfo North-Holland%Edited by
\publ Deutscher Verlag der Wiss. \vol \nofrills
\publaddr  Berlin
\yr  1978. 451 p.
\endref

\ref \no9 \by Reinov O.I.  \pages 43-47
\paper  Approximation properties of order p and the existence of
  non-p-nuclear operators with p-nuclear second adjoints
\yr 1981 \vol 256 \issue  1
\jour Doklady AN SSSR
\finalinfo
\endref

\ref \no10 \by Reinov O.I.\pages   125-134
\paper   Approximation properties of order p and the existence of
  non-p-nuclear operators with p-nuclear second adjoints
\yr 1982 \vol  109
\jour    Math. Nachr.
\endref

\ref \no11 \by Reinov O.I. \pages 145-165
\paper  Disappering tensor elements in the scale of p-nuclear operators
%\pages 145-165
\yr 1983\vol
\jour Theory of operators and theory of functions (LGU)
\endref


\ref \no12 \by Saphar P.\pages 71-100   $MR 43 \# 878.$
 \paper  Produits tensoriels d'espaces de Banach et classes
d'applications lineaires
 \yr 1970\vol  38
 \jour  Studia Math.
 \endref


\endRefs
\enddocument